\documentclass[12pt,a4paper]{amsart}
\usepackage{graphicx,amssymb}
\usepackage{amsmath,amssymb,amscd}
\input xy
\xyoption{all}
\usepackage[all]{xy}
\usepackage{hyperref}

\newcommand{\ra}{\rightarrow}

\newcommand{\ZZ}{\mathbb Z}

\newcommand{\PP}{\mathbb P}
\newcommand{\NN}{\mathbb N}

\newcommand{\cO}{\mathcal{O}}

\newcommand{\cM}{\mathcal{M}}

\newcommand{\Pic}{\mbox{Pic}}

\theoremstyle{plain}
\newtheorem{theorem}{Theorem}[section]
\newtheorem{lem}[theorem]{Lemma}
\newtheorem{prop}[theorem]{Proposition}
\newtheorem{cor}[theorem]{Corollary}

\theoremstyle{definition}
\newtheorem{rem}[theorem]{Remark}

\numberwithin{equation}{section}
\begin{document}
\title[Stable bundles with 4 sections]{Further examples of stable bundles of rank 2 with 4 sections}

\author{H. Lange}
\author{P. E. Newstead}

\address{H. Lange\\Department Mathematik\\
              Universit\"at Erlangen-N\"urnberg\\
              Bismarckstra\ss e $1\frac{ 1}{2}$\\
              D-$91054$ Erlangen\\
             Germany}
              \email{lange@mi.uni-erlangen.de}
\address{P.E. Newstead\\Department of Mathematical Sciences\\
              University of Liverpool\\
              Peach Street, Liverpool L69 7ZL, UK}
\email{newstead@liv.ac.uk}

\thanks{Both authors are members of the research group VBAC (Vector Bundles on Algebraic Curves). The second author 
would like to thank the Department Mathematik der Universit\"at 
         Erlangen-N\"urnberg for its hospitality}

\maketitle

\begin{center}
{\it Dedicated to the memory of Eckart Viehweg}
\end{center}

\begin{abstract}
In this paper we construct new examples of stable bundles of rank 2 of small degree with 4 sections on a smooth irreducible curve 
of maximal Clifford index. The corresponding Brill-Noether loci have negative expected 
dimension of arbitrarily large absolute value.  
\end{abstract}

\section{Introduction}

It has been apparent for some time that the classical Brill-Noether theory for line bundles on a smoooth 
irreducible curve does not extend readily to bundles of higher rank. Some aspects of this have been clarified 
recently by the introduction of Clifford indices of higher rank \cite{ln}. An example of a stable rank-3 bundle with 
Clifford index less than the classical Clifford index on a general curve of genus 9 or 11 is given in \cite{lmn}, 
disproving a conjecture of Mercat \cite{m}. Very recently, it was proved in \cite{fo} that there exist curves 
of any genus $\geq 11$ for which the rank-2 Clifford index is strictly smaller than the classical Clifford index. 
In this paper we use the methods 
of \cite{fo} to present further examples of this, showing in particular that the difference between the two Clifford 
indices can be arbitrarily large.\\

For any positive integer $n$ the rank-$n$ Clifford index $\gamma'_n(C)$ of a smooth projective curve of genus $g \geq 4$ 
over an algebraically closed field of characteristic 0 is defined as follows. For any vector bundle $E$ of rank $n$ 
and degree $d$ on $C$ define 
$$
\gamma(E) := \frac{1}{n}(d -2(h^0(E) -n)).
$$
Then 
$$
\gamma'_n = \gamma'_n(C) := \min \left\{\gamma(E) \;{\Big|} \begin{array}{c}
                                       E \; \mbox{semistable of rank} \; n \; \mbox{with} \\
                                        d \leq n(g-1) \; \mbox{and} \;  h^0(E) \geq 2n
                                        \end{array} \right\}.
$$   
Here $\gamma_1 = \gamma'_1$ is the classical Clifford index of $C$ and it is easy to see that $\gamma'_n \leq \gamma_1$ 
for all $n$. 

The gonality sequence $(d_r)_{r \in \NN}$ is defined by
$$
d_r
 := \min_{L \in Pic (C)} \{ \deg L \; | \; h^0(L) \geq r+1 \}.
$$
In classical terms $d_r$ is the minimum number $d$ for which a $g^r_d$ exists. In the case of a general curve we have for all $r$,
$$
d_r = g + r - \left[ \frac{g}{r+1} \right].
$$

According to \cite{m}, \cite{ln} a version of Mercat's conjecture states that
$$
\gamma'_n =\gamma_1 \quad \mbox{for all} \; n.
$$
As mentioned above, counterexamples in rank 3 and rank 2 are now known. For the rest of the paper we concentrate on rank 2. 

For $\gamma_1 \leq 4$ it is known that $\gamma'_2 = \gamma_1$ (see \cite[Proposition 3.8]{ln}). 
In any case, we have according to \cite[Theorem 5.2]{ln} 
$$
\gamma'_2 \geq \min \left\{ \gamma_1, \frac{d_4}{2} - 2 \right\}.
$$
For the general curve of genus 11 we have $\gamma_1 = 5$ and $d_4 = 13$. So in this case, $\gamma'_2 = 5$ or $\frac{9}{2}$. 
It is shown in \cite[Theorem 3.6]{fo} that there exist curves $C$ of genus 11 with $\gamma_1 = 5$ and $\gamma'_2 = \frac{9}{2}$,
but this cannot happen on a general curve of genus 11 \cite[Theorems 1.6 and 1.7]{fo}. Counterexamples to the conjecture in higher genus were also constructed 
in \cite{fo}. All examples $E$ constructed in \cite{fo} have $\gamma(E) = \gamma_1 - \frac{1}{2}$.

In this paper we use the methods of \cite{fo} to generalize these examples. Our main result is the following theorem.

\begin{theorem} \label{thm1.1}
Suppose $d = g-s$ with an integer $s \geq -1$ and 
$$
g \geq \max \{ 4s + 14,12\}.
$$
Suppose further that the quadratic form 
$$
3m^2 +dmn+ (g-1)n^2
$$
cannot take the value $-1$ for any integers $m,n \in \ZZ$. Then there exists a curve $C$ of genus $g$ having 
$\gamma_1(C) = \left[ \frac{g-1}{2} \right]$ and a stable bundle $E$ of rank $2$ on $C$ with 
$\gamma(E) = \frac{g- s}{2} -2$ and hence
$$
\gamma_1 - \gamma'_2 \geq \left[ \frac{g-1}{2} \right] - \frac{g - s}{2} + 2 > 0.
$$
In particular the difference $\gamma_1 - \gamma'_2$ can be arbitrarily large.
\end{theorem}

This statement can also be written in terms of the Brill-Noether loci $B(2,d,4)$ which are defined as follows. Let $M(2,d)$
denote the moduli space of stable bundles of rank 2 and degree $d$ on $C$. Then
$$
B(2,d,4) := \{ E \in M(2,d) \;|\; h^0(E) \geq 4 \}.
$$
Theorem \ref{thm1.1} says that under the given hypotheses $B(2,g-s,4)$ is non-empty. It may be noted that the expected 
dimension of $B(2,g-s,4)$ is $-4s-11 < 0$.

The key point in proving this theorem is the construction of the curves $C$, which all lie on K3-surfaces and are therefore
not general, although they do have maximal Clifford index.

\begin{theorem} \label{thm1.2}
Suppose $d = g-s$ with an integer $s \geq -1$ and 
$$
g \geq \max \{ 4s + 14,12\}.
$$
Then there exists a smooth K3-surface $S$ of type $(2,3)$ in $\PP^4$ containing a smooth curve $C$ of genus 
$g$ and degree $d$ with 
$$
\Pic (S) = H \ZZ \oplus C \ZZ,
$$
where $H$ is the polarization, such that $S$ contains no divisor $D$ with $D^2 = 0$.
Moreover, if $S$ does not contain a $(-2)$-curve, then $C$ is of maximal Clifford index $\left[ \frac{g-1}{2} \right]$.
\end{theorem}

The proof of Theorem \ref{thm1.2}, which uses the methods of \cite{f} and \cite{fo}, is given in Section 2. 
This is followed in Section 3
by the proof of Theorem \ref{thm1.1}.

\section{Proof of Theorem 1.2}

\begin{lem} \label{lem2.1}
Let $d = g-s$ with $g \geq 4s + 14$ and $s \geq -1$. Then $d^2 - 6(2g-2)$ is not a perfect square.
\end{lem}

\begin{proof}
If $d^2 -6(2g-2) = g^2 - (2s + 12)g + s^2 +12 = m^2$ for some non-negative integer $m$, then
the discriminant
$$
(s+6)^2 - (s^2 + 12 - m^2) = 12 s + 24 + m^2
$$
is a perfect square of the form $(m+ b)^2$ with $b \geq 2$. Solving the equation $g^2 - (2s + 12)g + (s^2 +12 - m^2) = 0$
for $g$, we get
\begin{equation} \label{eq2.1}
g = s+6 \pm (m + b).
\end{equation}
Now, since $b \geq 2$, we have $(m + b - 2)^2 \geq m^2$ and hence
$$
4(m+b) - 4 = (m+b)^2 - (m+b-2)^2 \leq 12 s + 24
$$
which gives $m+b \leq 3s + 7$. So \eqref{eq2.1} implies $g \leq 4s + 13$, which contradicts the hypothesis.
\end{proof}

\begin{prop} \label{prop2.2}
Let $g \geq 4s + 14$ with $s \geq -1$. Then there exists a smooth K3-surface $S$ of type $(2,3)$ in $\PP^4$ 
containing a smooth curve $C$ of genus $g$ and degree $d = g-s$ with 
$$
\Pic(S) = H \ZZ \oplus C \ZZ,
$$
where $H$ is the polarization, such that $S$ contains no divisor $D$ with $D^2 = 0$.
\end{prop}

\begin{proof}
The conditions of \cite[Theorem 6.1,2.]{k} are fulfilled to give the existence of $S$ and $C$. Let
$$
D = mH + nC \quad \mbox{with} \quad m,n \in \ZZ.
$$
We want to show that the equation $D^2 = 0$ does not have an integer solution. Now
$$
D^2 = 6m^2 + 2dmn + (2g-2)n^2.
$$
For an integer solution we must have that the discriminant $d^2 - 6(2g-2)$ is a perfect square and this contradicts Lemma
\ref{lem2.1}.
\end{proof}

\begin{lem} \label{lem2.3}
Under the hypotheses of Proposition 2.2, the curve $C$ is an ample divisor on $S$.
\end{lem}

\begin{proof}
We show that $C \cdot D > 0$ for any effective divisor on $S$ which we may assume to be irreducible.
So let $D \sim mH + nC$ be an irreducible curve on $S$.
So 
$$
C \cdot D = m(g-s) + n(2g-2).
$$ 
Note first that, since $H$ is a hyperplane, we have
\begin{equation} \label{eq2.2}
D \cdot H = 6m + (g-s)n  > 0.
\end{equation}
If $m,n \geq 0$, then one of them has to be positive and then clearly $C \cdot D > 0$. The case $m,n \leq 0$ 
contradicts \eqref{eq2.2}.

Suppose $m >0$ and $n < 0$.  Then, using \eqref{eq2.2} we have
$$
C \cdot D = m(g-s) + n(2g-2) > - n \left(\frac{(g-s)^2}{6} - (2g - 2) \right). 
$$ 
So $C \cdot D > 0$ for $g > s+6 + 2 \sqrt{3s + 6}$, which holds, since $g \geq 4s+14$.
 
Finally, suppose $m< 0$ and $n > 0$. Then, since we assumed $D$ irreducible,
$$
nC \cdot D = -m D \cdot H + D^2 \geq -m D \cdot H -2 \geq -m -2.
$$
If $m \leq -3$ , then $n C \cdot D > 0$.  
If $m = -1$, we have
$$
C \cdot D = -(g-s) + n(2g-2) \geq g +s -2 > 0. 
$$
The same argument works for $m = -2,\; n \geq 2$. Finally, if $m = -2$ and $n = 1$, we still get $C \cdot D > 0$
unless $D \cdot H = 1$ and $D^2 = -2$. Solving these equations gives $s = 1, g = 14$, contradicting the hypotheses.
\end{proof}

\begin{theorem} \label{thm2.4}
Let the situation be as above with $d = g-s$, $s \geq -1$ and 
$$
g \geq \max \{ 4s +14,12 \}.
$$
If $S$ does not contain a $(-2)$-curve, then $C$ is of maximal Clifford index $\left[ \frac{g-1}{2} \right]$.
\end{theorem}

Note that a stronger form of this has been proved for $s = -2$ and $g$ odd in \cite[Theorem 3.6]{fo} and for $s = -1$ 
and $g$ even in \cite[Theorem 3.7]{fo}. The proof follows closely that of \cite[Theorem 3.3]{f}, but, since some of the estimates 
are delicate and our hypotheses differ, we give full details.

\begin{proof}
Since $C$ is ample by Lemma \ref{lem2.3}, it follows from \cite[Proposition 3.3]{cp} that $C$ is of Clifford dimension 1.

Suppose that $\gamma_1(C) < \left[\frac{g-1}{2} \right]$. According to \cite{dm} there is an effective divisor $D$ on $S$ 
such that $D|_C$ computes $\gamma_1(C)$ and satisfying
$$
h^0(S,D) \geq 2, \quad h^0(S,C-D) \geq 2 \quad  \mbox{and} \quad \deg(D|_C) \leq g-1.
$$ 
We consider the exact cohomology sequence
$$
0 \ra H^0(S,D-C) \ra H^0(S,D) \ra H^0(C,D|_C) \ra H^0(S,D-C).
$$
Since $C -D$ is effective, and not equivalent to zero, we get
$$
H^0(S,D-C) = 0.
$$
By assumption $S$ does not contain $(-2)$-curves, so $|D-C|$ has no fixed components. 
According to Proposition \ref{prop2.2} the equation $(C -D)^2 = 0$ has no solutions, therefore
$(C-D)^2 >0$ and the general element of $|C-D|$ is smooth and irreducible. It follows that 
$$
H^1(S,D-C) = H^1(S,C-D)^* = 0
$$
and 
$$
\gamma_1(C) = \gamma(D|_C) = D \cdot C - 2 \dim |D| = D \cdot C - D^2 -2
$$
by Riemann-Roch. We shall prove that 
$$
D \cdot C - D^2 - 2 \geq \left[ \frac{g-1}{2} \right],
$$
a contradiction.\\

Let $D \sim mH + nC$ with $m,n \in \ZZ$. Since $D$ is effective and $S$ contains no $(-2)$-curves, 
we have $D^2 > 0$ and $D \cdot H > 2$. Since $C-D$ is also effective, we have $(C-D)\cdot H > 2$, i.e. 
$D \cdot H < d - 2$. These inequalities and $\deg(D|_C) \leq g-1$ translate to the following inequalities
\begin{equation} \label{eq3.1}
3m^2 + mnd + n^2(g-1) > 0,
\end{equation}  
\begin{equation} \label{eq3.2}
2 < 6m +nd < d-2,
\end{equation}
\begin{equation} \label{eq3.3}
md + (2n-1)(g-1) \leq 0.
\end{equation}
Consider the function
$$
f(m,n) := D\cdot C -D^2 -2 = -6m^2 + (1 -2n)dm + (n-n^2)(2g-2) -2,
$$
and denote by
$$
a := \frac{1}{6}(d + \sqrt{d^2 -12(g-1)}) \quad \mbox{and}  \quad b:= \frac{1}{6}(d - \sqrt{d^2 -12(g-1)})
$$
the solutions of the equation $6x^2 - 2dx + 2g-2 = 0$. 
Note that $d^2 > 12(g-1)$. So $a$ and $b$ are positive real numbers.

Suppose first that $n < 0$. From \eqref{eq3.1} we have either $m < -bn$ or $m > -an$. If $m < -bn$, then \eqref{eq3.2}
implies that $2 < n(d-6b) < 0$, because $n < 0$ and $d -6b = \sqrt{d^2 - 12(g-1)} > 0$, which gives a contradiction.

If $n < 0$ and $m > -an$, from \eqref{eq3.3} we get
$$
-an < m \leq \frac{(g-1)(1-2n)}{d} < \frac{(1-2n)d}{12},
$$
since $d^2 > 12(g-1)$. For a fixed $n$, $f(m,n)$ is increasing as a function of $m$ for $m \leq \frac{(1-2n)d}{12}$
and therefore
\begin{eqnarray*}
f(m,n) &>& f(-an,n)\\
& = & \frac{d^2 - 12(g-1) +d\sqrt{d^2 - 12(g-1)}}{6} \cdot(-n) -2 \\
& \geq & \frac{d^2 - 12(g-1) +d\sqrt{d^2 - 12(g-1)}}{6}  -2 \\
&\geq& \frac{g-1}{2},
\end{eqnarray*}
which gives a contradiction. Here the last inequality reduces to 
$$
d \sqrt{d^2 -12(g-1)} \geq 15g -3 - d^2
$$ 
which certainly holds if $d^2 \geq 15g -3$. This is true under our hypotheses on $g$ if $s \geq 1$. The inequality can be 
checked directly in the cases $s = 0$ and $s = -1$. 

Now suppose $n > 0$. From \eqref{eq3.1} we get that either $m < -an$ or $m > -bn$. If $m < -an$, we get from \eqref{eq3.2},
$2 < n(-6a +d) < 0$, a contradiction.

When $m > -bn$, first suppose $n = 1$.
Then \eqref{eq3.3} gives 
\begin{equation} \label{eq3.4}
-b < m \leq -\frac{g-1}{d}.
\end{equation}
We claim that
\begin{equation} \label{eq3.5}
1 < b < \frac{4}{3}.
\end{equation}
In terms of $s$ we have
\begin{eqnarray*}
6b & = & g - s - \sqrt{(g- s)^2 - 12(g-1)}\\
& = & g-s - \sqrt{(g -(s + 6))^2 - 12 s -24}\\
& > & g - s - (g -(s + 6)) = 6,
\end{eqnarray*}
since $s \geq -1$. This gives $1 < b$. For the second inequality note that $b = \frac{4}{3}$ gives 
$s = \frac{g-13}{4}$ and $b$ is a strictly increasing function of $s$ in the interval 
$\left[ -1,  \frac{g-13}{4} \right]$. Since certainly $s < \frac{g-13}{4}$, we obtain
$b < \frac{4}{3}$. 

So there are no solutions of \eqref{eq3.4} unless $d \geq g-1$, i.e. $s = 1,0$ or $-1$.
For these values of $s$ we must have $m = -1$ and 
$$
f(m,n) = f(-1,1) = d-8.
$$
So $f(-1,1) \geq \left[ \frac{g-1}{2} \right]$ if and only if $g \geq 2s + 14$.

Now suppose $m > -bn$ and $n \geq 2$. Then \eqref{eq3.3} gives
$$
f(m,n) \geq \min \left\{ f \left( -\frac{(g-1)(2n-1)}{d},n \right), f(-bn,n) \right\}.
$$ 
We have
$$
f \left( -\frac{(g-1)(2n-1)}{d},n \right) = \frac{g-1}{2} \left((2n-1)^2 \left( 1 - \frac{12(g-1)}{d^2} \right) +1 \right) -2.
$$
It is easy to see that $f \left( -\frac{(g-1)(2n-1)}{d},n \right) \geq \frac{g-1}{2}$ for $n \geq 2$.
Moreover,
$$
f(-bn,n) = -bdn + n(2g-2)-2 = n(2g-2 -bd) -2.
$$
Note that
$$
2g-2-bd = \frac{\sqrt{d^2-12(g-1)}}{6}(d - \sqrt{d^2-12(g-1)}) > 0.
$$
So $f(-bn,n)$ is a strictly increasing function of $n$. Hence it suffices to show that
$f(-2b,2) \geq \frac{g-1}{2}$ or equivalently
$$
7(g-1) - 4bd - 4 \geq 0.
$$
According to \eqref{eq3.5} we have $b < \frac{4}{3}$. So, since $d \leq g+1$, we have
\begin{eqnarray*}
7(g-1) - 4bd -4 & \geq & 7(g-1) - \frac{16}{3}d -4 \\
& \geq & 7g-7 - \frac{16}{3}g - \frac{16}{3} -4 = \frac{1}{3}(5g-49) > 0.
\end{eqnarray*}
This completes the argument for $m > -bn, \; n > 0$.

Finally, suppose $n=0$. Then
$$
f(m,0) = -6m^2 + dm -2.
$$
As a function of $m$ this takes its maximum value at $\frac{d}{12}$. By \eqref{eq3.3}, 
$m \leq \frac{g-1}{d} \leq \frac{d}{12}$. 
So $f(m,0)$ takes its minimal value in the allowable range at $m=1$. Since $f(1,0) = d-8$, we require 
$d-8 \geq \left[ \frac{g-1}{2}  \right]$
or equivalently
$$
g \geq 2 s + 14,
$$
which is valid by hypothesis.
\end{proof}

This completes the proof of Theorem \ref{thm1.2}.

\begin{rem} \label{rem2.5}
For $s = 0$ or $-1$ the assumptions of the theorem are best possible, since in these cases 
$\gamma(H|_C) = \gamma((C-H)|_C) = d-8$ would otherwise 
be less than $\left[ \frac{g-1}{2} \right]$. For $s \geq 1$ the conditions can be relaxed. 
For example, if $s \geq 1$ and $g = 4s + 12$, the only places where the argument can fail
are in the proofs of Lemma \ref{lem2.1} and formula \eqref{eq3.5}. In the first case, one can show directly 
that $d^2 -6(2g-2)$ is not a perfect square; in the second, one can show that $b < \frac{3}{2}$, which is sufficient. 
\end{rem}

\begin{rem} \label{rem2.6}
The condition that $S$ does not contain a $(-2)$-curve certainly holds if $3m^2 +dmn + (g-1)n^2 = -1$ has no
solutions. We do not know precisely when this is true, but it certainly holds if both $g-1$ and $g - s$ are divisible by 3.
So the conclusion of Theorem \ref{thm2.4} holds for $s \equiv 1 \mod 3$, if $g \geq 4s + 14$ and $g \equiv 1 \mod 3$. 
The conclusion also holds, for example, for $g=16$ and $s = 1$ (see Remark \ref{rem2.5}). 
\end{rem}

\section{Proof of Theorem 1.1}

\begin{lem} \label{lem3.1}
Let $C$ and $H$ be as in Proposition \ref{prop2.2} with $d = g- s, s \geq -1$ and suppose that $S$ has no $(-2)$-curves.
Then $H|_C$ is a generated line bundle on $C$ with $h^0(\cO_C(H|_C)) = 5$ and 
$$
S^2H^0(\cO_C(H|_C)) \ra H^0(\cO_C(H^2|_C))
$$ 
is not injective.
\end{lem}

\begin{proof}
Consider the exact sequence
$$
0 \ra \cO_S(H-C) \ra \cO_S(H) \ra \cO_C(H|_C) \ra 0.
$$
$H -C$ is not effective, since $(H-C) \cdot H = 6 -d < 0$. So we have 
$$
0 \ra H^0(\cO_S(H)) \ra H^0(\cO_C(H|_C)) \ra H^1(\cO_S(H-C)) \ra 0.
$$
Now 
$$
(C-H)^2 = 2g-2 -2d + 6 = 2 s + 4 \geq 2
$$
and 
$$
H^2(\cO_S(C-H)) = H^0(\cO_S(H-C))^* = 0.
$$
So by Riemann-Roch $h^0(\cO_S(C-H)) \geq 3$. Since $S$ has no $(-2)$-curves, it follows that the linear system $|C-H|$
has no fixed components and hence its general element is smooth and irreducible (see \cite{sd}). Hence $h^1(\cO_S(H-C)) = 0$ and 
therefore
$h^0(\cO_C(H|_C)) = h^0(\cO_S(H)) =5$. The last assertion follows from the fact that $S$ is contained in a quadric. 
\end{proof}

\begin{rem}
Lemma \ref{lem3.1} implies that $H|_C$ belongs to $W^4_{g- s}$. So $g- s \geq d_4$. Since the generic 
value of $d_4$ is $g +4 - \left[ \frac{g}{5} \right]$, it follows that $C$ has non-generic $d_4$ if $g < 5 s + 20$.
\end{rem}

\begin{lem} \label{lem3.3}
Let $C$ be a smooth irreducible curve and $M$ a generated line bundle on $C$ of degree $d < 2d_1$ with $h^0(M) = 5$
and such that $S^2H^0(M) \ra H^0(M^2)$ is not injective.
Then $B(2,d,4) \neq \emptyset$.
\end{lem}

The proof is identical with that of \cite[Theorem 3.2 (ii)]{gmn}. $\hspace{2.1cm} \square$

\begin{theorem} \label{thm3.4}
Let $C$ be as in Theorem \ref{thm2.4}. Then
\begin{itemize}
\item[(i)] $B(2,g-s,4) \neq \emptyset$;
\item[(ii)] $\gamma'_2(C) \leq \frac{g- s}{2} - 2 < \gamma_1(C)$.
\end{itemize}
\end{theorem}

\begin{proof}
This follows from Theorem \ref{thm2.4} and Lemmas \ref{lem3.1} and \ref{lem3.3}.
\end{proof}

This completes the proof of Theorem \ref{thm1.1}, where the last assertion follows from Remark \ref{rem2.6}. 

\begin{cor}
$\gamma'_{2n}(C) < \gamma_1(C)$ for every positive integer $n$.
\end{cor}

\begin{proof}
This follows from Theorem \ref{thm3.4} and \cite[Lemma 2.2]{ln}.
\end{proof}

\begin{rem}
Under the conditions of Theorem \ref{thm1.1}, for any stable bundle $E$ of rank 2 and degree $g-s$ on $C$ with
$h^0(E) = 4$, it follows from \cite[Proposition 5.1]{gmn} that the coherent system $(E,H^0(E))$ is $\alpha$-stable
for all $\alpha > 0$. So the corresponding moduli spaces of coherent systems are non-empty.
\end{rem}

\end{document}